\documentclass[english,11pt]{article}
\usepackage{amsopn,amsthm,amsfonts,amsmath,amssymb,amscd}
\usepackage{babel}
\usepackage{amssymb}
\usepackage{graphicx}

\theoremstyle{definition}
\newtheorem{Theorem}{Theorem}
\newtheorem{Lemma}[Theorem]{Lemma}

\newtheorem{Corollary}[Theorem]{Corollary}
\newtheorem{Example}[Theorem]{Example}

\newtheorem*{Abstract}{Abstract}

\newenvironment{Proof*}{{\it Proof.}}

\newcommand{\NN}{\mathbb{N}}

\newcommand{\DEF}[1]{\emph{#1}}

\newcommand{\adots}{\mathinner{\mskip1mu\raise1pt\vbox{\kern7pt\hbox{.}}\mskip2mu\raise4pt\hbox{.}\mskip2mu
\raise7pt\hbox{.}\mskip1mu}}
\newcommand{\im}{\mathop{\rm im}}

\newcommand{\rk}{\mathop{\rm rk}}
\newcommand{\sh}[1]{\mathop{\rm sh}(#1)}
\newcommand{\lin}{{\mathcal L}in}

\newcommand{\n}{{\mathcal N}}
\newcommand{\nn}{{\mathcal N}_2}
\newcommand{\nschab}{{\mathcal V}(n,a,b)}
\newcommand{\nsch}{{\mathcal V}(n)}
\newcommand{\nbsch}{{\mathcal V}_{B}(n)}

\newcommand{\partition}[1]{{\mathcal P}(#1)}

\newcommand{\ord}[1]{\mathop{\rm ord}(#1)}

\newcommand{\bkp}[2]{\mathcal{B}(#1;#2)}
\newcommand{\sch}[2]{\mathcal{S}(#1;#2)}

\newcommand{\gl}[2]{\mathop{GL}_{#1}(#2)}

\begin{document}

\title{Jordan forms for mutually annihilating nilpotent pairs}
\author{Polona Oblak \\
 \small{Institute of Mathematics, Physics and Mechanics, Department of Mathematics,} \\
 \small{Jadranska 19, SI-1000 Ljubljana, Slovenia}
 \vspace{2mm}\\
 \small{polona.oblak@fmf.uni-lj.si}
}

\date{}

\maketitle 
\parindent=0cm

\Abstract

In this paper we completely characterize all possible pairs of Jordan canonical forms
for mutually annihilating nilpotent pairs, i.e. pairs $(A,B)$ of nilpotent matrices such that $AB=BA=0$.
\bigskip
\bigskip

{\small

{\bf Math. Subj. Class (2000):} 15A03, 15A21, 15A27

}

\vspace{5mm}

\bigskip
\bigskip

\section{Introduction} \label{sec:intro}

We consider pairs of $n \times n$ commuting matrices over an algebraically closed field $F$.
For $n,a,b$ (all at least 2) let $\nschab$ be the variety of all pairs $(A,B)$ of commuting nilpotent matrices such that $AB=BA=A^a=B^b=0$.
In \cite{schroeer} Schr\" oer classified the irreducible components of $\nschab$ and thus answered a question stated by
Kraft \cite[p. 201]{kraft} (see also \cite{gp} and \cite{laubenbacher}).

If $\underline{\mu}=(\mu_1,\mu_2,\ldots,\mu_t)$ is a partition of $n$ then we denote by ${\cal O}_{\underline{\mu}}$ the conjugacy class
of all nilpotent matrices such that the sizes of Jordan blocks in its Jordan canonical form are equal to $\mu_1, \mu_2,\ldots,\mu_t$.
Let $\underline{\mu}=(\mu_1,\mu_2,\ldots,\mu_t)$ and $\underline{\nu}=(\nu_1,\nu_2,\ldots,\nu_s)$ be partitions of $n$
such that $\mu_1 \leq a$ and $\nu_1 \leq b$ and let $\pi_1(A,B)=A$ and $\pi_2(A,B)=B$ for $(A,B) \in \nschab$ be the projection maps.
Schr\" oer \cite[p. 398]{schroeer} noted that the intersection of the fibers 
$\pi_1^{-1}({\cal O}_{\underline{\mu}}) \cap \pi_2^{-1}({\cal O}_{\underline{\nu}})$
is not very well-behaved for different reasons. It might be empty or reducible and the closure
of the intersection of fibers is in general not a union of such intersections.

In this paper we answer the question for which pairs of partitions $(\underline{\mu},\underline{\nu})$ the intersection 
$\pi_1^{-1}({\cal O}_{\underline{\mu}}) \cap \pi_2^{-1}({\cal O}_{\underline{\nu}})$ is nonempty. If 
$\underline{\mu}=(\mu_1,\mu_2,\ldots,\mu_t,1^m)$, where $\mu_t \geq 2$, is a fixed partition, then partitions $\underline{\nu}$, such that
$\pi_1^{-1}({\cal O}_{\underline{\mu}}) \cap \pi_2^{-1}({\cal O}_{\underline{\nu}})$ is nonempty, are of the form 
$$(\lambda_1+\varepsilon_1,\lambda_2+\varepsilon_2, \ldots,\lambda_l+\varepsilon_l,2^c,1^d),$$
where $(\lambda_1,\lambda_2,\ldots,\lambda_l)$ is a partition of $m$ and $\varepsilon_i \in \{0,1,2\}$. The precise constraints for $\varepsilon_i$, $c$ and $d$ are given in Theorem \ref{thm:nbsch}.

\medskip

The structure of varieties of commuting pairs of matrices and of commuting pairs of nilpotent matrices is not yet well understood. 
A motivation for our study is also to contribute to better understanding of the structure of the latter variety
and to help in understanding the (ir)reducibility of the variety of triples of commuting matrices 
(see also \cite{gerstenhaber}, \cite{guralnick}, \cite{gs}, \cite{han}, \cite{ho}, \cite{mt}, \cite{omladic}). 
It was proved by Motzkin and Taussky \cite{mt} (see also Gerstenhaber \cite{gerstenhaber} and Guralnick \cite{guralnick}), 
that the variety of pairs of commuting matrices was irreducible. Guralnick \cite{guralnick} was the first to show that this is no 
longer the case for the variety of triples of commuting matrices. Recently, it was proved that the variety of 
commuting pairs of nilpotent matrices was irreducible (Baranovsky \cite{baranovsky}, Basili \cite{basili}). Some of our 
results on Jordan canonical forms for commuting nilpotent pairs appear in \cite{oblak}.

\smallskip

Let us briefly describe the setup of the paper. The second section is of preparatory nature.
In Section 3, we consider ${\cal V}(n;n,n)$ and we fix the partition $\underline{\mu}$ corresponding to $B$. Then we describe a 
sequence of conjugations of matrix $A$, which puts $A$ to a nicer pattern, called the reduced form of $A$. For $A$ in
reduced form, we compute its Jordan canonical form in Section 4. As a corollary, we give for general $\nschab$ 
all pairs of partitions $({\underline{\mu}}, {\underline{\nu}})$ such that the intersection 
$\pi_1^{-1}({\cal O}_{\underline{\mu}}) \cap \pi_2^{-1}({\cal O}_{\underline{\nu}})$ is nonempty.

\bigskip
\bigskip

\section{Notation and first observations} 

\bigskip

Let us denote by $\n=\n(n,F)$ the variety of all $n \times n$ nilpotent matrices over 
an algebraically closed field $F$ and write $\nn=\{(A,B) \in \n \times \n; AB=BA\}$. 

We call a pair $(A,B) \in \nn$ a \DEF{mutually annihilating nilpotent pair} if $AB=BA=0$ and we denote the set 
of all mutually annihilating nilpotent pairs by $\nsch=\{(A,B) \in \nn; \; AB=BA=0\}$.
Additionally, we write $$\nbsch=\{A \in \n; \; (A,B) \in \nsch\}$$ for $B \in \n$.



\bigskip

Suppose that $\mu_1 \geq \mu_2 \geq \ldots \geq \mu_t > 0$ are
the sizes of Jordan blocks in the Jordan canonical form for $B$. 
We call the partition $\underline{\mu}=(\mu_1,\mu_2,\ldots,\mu_t)$ 
the \DEF{shape} of the matrix $B$ and denote it by $\sh{B}$. 
We also write $\sh{B}=(\mu_1,\mu_2,\ldots,\mu_t)=(m_1^{r_1},m_2^{r_2},\ldots,m_l^{r_l})$, where $r_i >0$ and $\sum_{i=1}^l r_i=t$. 
For a partition $\underline{\mu}=(\mu_1,\mu_2,\ldots,\mu_t)$ we write $|\underline{\mu}|=\sum_{i=1}^{t} \mu_i$ and
for a sequence $(a_1,a_2,\ldots,a_k)$, $a_i \in \NN$, we write
$\ord{a_1,a_2,\ldots,a_k}=(a_{\pi(1)},a_{\pi(2)},\ldots,a_{\pi(k)})$,  where  
$a_{\pi(1)}\geq a_{\pi(2)}\geq \ldots \geq a_{\pi(k)}$ and $\pi$ is a permutation of $\{1,2,\ldots,k\}$.

By $J_{\underline{\mu}}=J_{(\mu_1,\mu_2,\ldots,\mu_t)}=J_{\mu_1}\oplus J_{\mu_2}\oplus\ldots\oplus J_{\mu_t}$ we denote the 
upper triangular matrix in its Jordan canonical form with blocks of sizes $\mu_1 \geq \mu_2 \geq  \ldots \geq  \mu_t$.

\medskip

Let $\partition{n}$ denote the set of all partitions of $n \in \NN$
and for subsets ${\cal S}_1 \subseteq \n$ and ${\cal S}_2 \subseteq \nn$ we write 
$$\partition{{\cal S}_1}=\{\sh{A}; \; A\in {\cal S}_1 \} \; \text{ and } \;
\partition{{\cal S}_2}=\{(\sh{A},\sh{B}); \; (A,B) \in {\cal S}_2 \}\, .$$
Note that $(\underline{\nu},\underline{\mu}) \in \partition{\nn}$ (resp. in $\partition{\nsch}$) if and only if
$(\underline{\mu},\underline{\nu}) \in \partition{\nn}$ (resp. in $\partition{\nsch}$), i.e. $\partition{\nn}$ and 
$\partition{\nsch}$ are symmetric.

\bigskip

Let $(A,B) \in \nsch$ and $\sh{(A,B)}=(\underline{\nu},\underline{\mu})$. 
Then there exists $P \in \gl{n}{F}$ such that $\sh{(PAP^{-1},J_{\underline{\mu}})}=(\underline{\nu},\underline{\mu})$
and $(PAP^{-1},J_{\underline{\mu}})\in \nsch$.
Therefore we can assume that $B$ is already in its upper triangular Jordan canonical form.

\medskip

Suppose that $B=J_{(\mu_1,\mu_2,\ldots,\mu_t)}$ and write $A=[A_{ij}]$ where $A_{ij} \in {\mathcal M}_{\mu_i \times \mu_j}$. 
It is well known that if $AB=BA$, then $A_{ij}$ are all upper triangular Toeplitz matrices, i.e.
for $1 \leq i \leq j \leq t$ we have
\begin{equation}\label{eq:com}
A_{ij}=
\left[ \begin{matrix}
0 & \ldots & 0 & a_{ij}^0 & a_{ij}^1 & \ldots & a_{ij}^{\mu_i-1}\\
\vdots &  & \ddots & 0 & a_{ij}^0 & \ddots &  \vdots \\
\vdots &  &  & \ddots & 0 & \ddots &  a_{ij}^{1} \\
0 & \ldots & \ldots  & \ldots  & \ldots  & 0 & a_{ij}^{0}
\end{matrix}
\right]
\;
\text{ and }
\;
A_{ji}=
\left[ \begin{matrix}
a_{ji}^0 & a_{ji}^1 & \ldots & a_{ji}^{\mu_i-1}\\
0 & a_{ji}^0 & \ddots &  \vdots \\
\vdots &  \ddots & \ddots &  a_{ji}^{1} \\
\vdots &   & 0 &  a_{ji}^{0} \\
\vdots &   &  &  0 \\
0 & \ldots & \ldots   & 0
\end{matrix}
\right]
\, .
\end{equation}
If $\mu_{i}=\mu_j$ then we omit the rows or columns of zeros in $A_{ij}$ or $A_{ji}$ above.

\bigskip

\Lemma \label{thm:schblocks}
 If $(A,J_{\underline{\mu}}) \in \nsch$, where $\underline{\mu}=(\mu_1,\mu_2,\ldots\mu_t)$, then for all $1 \leq i, j \leq t$ we have
 $$A_{ij}=
 \left[ \begin{matrix}
 0 & \ldots & 0 &  a_{ij}^{\mu_i-1}\\
 \vdots &  & \ddots & 0 \\
 \vdots &  &   & \vdots \\
 0 & \ldots & \ldots  & 0
 \end{matrix}
 \right]
 $$
 for all $i,j$, such that $\mu_i, \mu_j \geq 2$.

 \medskip
 
 \begin{proof}
  Since $A$ and $J_{\underline{\mu}}$ commute, $A=[A_{ij}]$ is of the form \eqref{eq:com}.
  Since they are mutually annihilating, we have $A_{ij}J_{\mu_j}=0$ and $J_{\mu_i}A_{ij}=0$ for all $1 \leq i, j \leq t$
  and thus $a_{ij}^s=0$ for all $0 \leq s \leq \min\{\mu_i,\mu_j\}-2$.   
 \end{proof}
 
\bigskip

Let us write $\underline{\mu}=(\mu_1,\mu_2,\ldots,\mu_k,1^m)$, where $\mu_k \geq 2$.

\Theorem\label{thm:1}
Let $(A,B)$ be a mutually annihilating nilpotent pair.
\begin{enumerate}
 \item If $\sh{B}=(1^n)$, i.e. $B=0$, then $\partition{\nbsch}=\partition{n}$.
 \item If $\sh{B}=(\mu_1,\mu_2,\ldots,\mu_k)$, where $\mu_k \geq 2$, then $\partition{\nbsch}=\{(2^i,1^{n-2i}); \; i=0,1,\ldots,k\}$. 
\end{enumerate}

\begin{Proof*} 
\begin{enumerate}
 \item If $\sh{B}=(1^n)$, then $B=0$ and $(A,B)$ is a mutually annihilating nilpotent pair for every nilpotent matrix $A$. Thus 
       $\partition{\nbsch}=\partition{\n}=\partition{n}$.
 \item Suppose that $B=J_{\underline{\mu}}$, where $\underline{\mu}=(\mu_1,\mu_2,\ldots,\mu_k)$ and $\mu_k \geq 2$. Then 
       $A$ has a form of Lemma \ref{thm:schblocks}, with $\mu_i\geq 2$ for all $i$ 
       and thus $A^2=0$. 
       Since $\rk\left(A\right)$ can be 0, 1, ..., or $k$, it follows that $\sh{A}=(2^i,1^{n-2i})$ for some 
       $0 \leq i \leq k$. Therefore $\partition{\nbsch}=\{(2^i,1^{n-2i}); \; i=0,1,\ldots,k\}$. 
       \hfill$\blacksquare$
\end{enumerate} 
\end{Proof*}
  
\bigskip
\bigskip

\section{Reduced form for $A$ under similarity} 

\bigskip

We denote by $\bkp{\underline{\mu}}{\underline{\nu}}$ the set of all $|\underline{\mu}|\times|\underline{\nu}|$ 
matrices, which are block upper triangular Toeplitz when partitioned according to row partition 
$\underline{\mu}$ and column partition $\underline{\nu}$.
By $\sch{\underline{\mu}}{\underline{\nu}}$ we denote all the matrices from $\bkp{\underline{\mu}}{\underline{\nu}}$
with at most one nonzero entry in each block.

For example, $A \in \bkp{\underline{\mu}}{\underline{\mu}}$ if and only if
$A$ commutes with $J_{\underline{\mu}}$ and $A \in \sch{\underline{\mu}}{\underline{\mu}}$ if and only if
$(A,J_{\underline{\mu}})$ is a mutually annihilating nilpotent pair (see Lemma \ref{thm:schblocks}).

\bigskip

In order to find a matrix that is similar to $A$ and has nicer pattern, let us consider the conjugate action of $\gl{n}{F}$ on 
$\sch{\underline{\mu}}{\underline{\mu}}$. 


\bigskip

Fix a matrix $B=J_{(\underline{\mu},1^m)}$, with $\underline{\mu}=(\mu_1,\mu_2,\ldots,\mu_k)$, where 
$\mu_k \geq 2$ and $m \geq 1$, and let $(A,B)$ be a mutually annihilating nilpotent pair. Then $A$ can be written as 
$A=\left[ \begin{matrix}
 A_{11} & A_{12}\\
 A_{21} & A_{22}
\end{matrix} \right]$, where $A_{11} \in \sch{\underline{\mu}}{\underline{\mu}}$, 
$A_{12} \in \sch{\underline{\mu}}{1^m}$, $A_{21} \in \sch{1^m}{\underline{\mu}}$ and
$A_{22} \in \sch{1^m}{1^m}$. According to Basilli  \cite[Lemma 2.3]{basili}, $A_{22}$ is nilpotent and thus there exists $P \in \gl{m}{F}$ 
such that $P^{-1}A_{22}P=J_{(\lambda_1,\lambda_2,\ldots,\lambda_l)}=J_{\underline{\lambda}}$.
Then 
$$\left[ \begin{matrix}
 I & 0\\
 0 & P^{-1}
\end{matrix} \right]
\,
\left[ \begin{matrix}
 A_{11} & A_{12}\\
 A_{21} & A_{22}
\end{matrix} \right]
\,
\left[ \begin{matrix}
 I & 0\\
 0 & P
\end{matrix} \right]=
\left[ \begin{matrix}
 A_{11} & A_{12}P\\
 P^{-1}A_{21} & J_{\underline{\lambda}}
\end{matrix} \right] \, .$$
Observe that $A_{12}P \in \sch{\underline{\mu}}{1^m}$ and $P^{-1}A_{21} \in \sch{1^m}{\underline{\mu}}$.
To describe $\partition{\nbsch}$ it suffices to describe all possible partitions for nilpotent matrices $\left[ \begin{matrix}
 A_{11} & A_{12}P\\
 P^{-1}A_{21} & J_{\underline{\lambda}}
\end{matrix} \right]
$.
From now on, we assume that $A$ is already in this form.

\bigskip

We denote by $E_{i,r_i;j,r_j}(\xi)$ the $n \times n$ matrix, partitioned according to $\underline{\nu} \times \underline{\nu}$, which 
has the diagonal entries equal to 1 and the entry in the $(\nu_1+\nu_2+\ldots+\nu_{i-1}+r_i)$-th row and $(\nu_1+\nu_2+\ldots+\nu_{j-1}+r_j)$-th 
column equals to $\xi$. All other entries of $E_{i,r_i;j,r_j}(\xi)$ are equal to 0.

The conjugation of $A$ by the matrix $E_{i,r_i;j,r_j}(\xi)$ corresponds to the following operations on rows and columns of $A$:
we add the $r_i^{\rm th}$ row in the $i^{\rm th}$ block of rows, multiplied by $\xi$, to the $r_j^{\rm th}$th row in the $j^{\rm th}$ 
block of rows of matrix $A$ and at the same time we add the $r_j^{\rm th}$ column of the $j^{\rm th}$ block of columns, multiplied by $-\xi$, 
to the $r_i^{\rm th}$ column in the $i^{\rm th}$ block of columns of $A$.

\bigskip

Now, let us consider some special conjugations of $A$ by matrices $E_{i,r_i;j,r_j}(\xi)$.

\begin{enumerate}
 \item Conjugation by $E_{k+j,i,t,1}(\xi)$, where $1 \leq j \leq l$, $1 \leq i \leq \lambda_j-1$, $1 \leq t \leq k$,
       adds a row with at most $k+1$ nonzero elements (and and least one 1) to the first row in some block of rows. Since all the entries
       of the corresponding first column of a block of columns are equal to 0, the conjugation by $E_{k+j,i,t,1}(\xi)$ 
       affects the matrix $A$ only by its row operation. Thus $A_{21}$ and $A_{22}$ do not change. 
       The entries of $A_{11} \in \bkp{\underline{\mu}}{\underline{\mu}}$ may change, but $A_{11}$ remains in 
       $\sch{\underline{\mu}}{\underline{\mu}}$ after the conjugation.
       
       Hence, it is possible to choose $\xi$ such that the conjugated matrix $A$ has 
       $A(\mu_1+\mu_2+\ldots+\mu_{t-1}+1;n-m+\lambda_1+\lambda_2+\ldots+\lambda_{j-1}+i)=0$
       for all $1 \leq j \leq l$, $1 \leq t \leq k$ and $2 \leq i \leq \lambda_j$. It means that we may assume that the only nonzero entries
       of $A_{12} \in \bkp{\underline{\mu}}{\underline{\lambda}}$ are the ones in the top left corner of each block according to partition 
       $(\underline{\mu}, \underline{\lambda})$.

 \item Similarly as in 1., conjugation of $A$ by $E_{t,\mu_t,k+j,i}(\xi)$, where $1 \leq t \leq k$, $1 \leq j \leq l$ and $2 \leq i \leq \lambda_j$
       affects matrix $A$ only by its column operation and thus $A_{11}$ (note that $i \geq 2$), $A_{12}$ and $A_{22}$ do not change.
      
       Hence, it is possible to choose $\xi$ such that the conjugated matrix $A$ has 
       $A(n-m+\lambda_1+\lambda_2+\ldots+\lambda_{j-1}+i;\mu_1+\mu_2+\ldots+\mu_{t})=0$
       for all $1 \leq j \leq l$, $1 \leq t \leq k$ and $1 \leq i \leq \lambda_j-1$. It means that we may assume that the only nonzero entries
       of $A_{21} \in \bkp{\underline{\lambda}}{\underline{\mu}}$ are the ones in the lower right corner of each block according
       to partition $(\underline{\lambda}, \underline{\mu})$.
\end{enumerate}

From now on, we assume that the only possible nonzero entries of block $A_{12} \in \bkp{\underline{\mu}}{\underline{\lambda}}$ of matrix $A$
are the ones in the top left corner of each block according to partition $(\underline{\mu}, \underline{\lambda})$
and the only possible nonzero entries
of $A_{21} \in \bkp{\underline{\lambda}}{\underline{\mu}}$ are the ones in the lower right corner of each block according
to partition $(\underline{\lambda}, \underline{\mu})$.

\begin{enumerate}
 \item[3.] Suppose that there exist $1 \leq t \leq k$ and $1 \leq s \leq l$ such that 
       $A(\mu_1+\mu_2+\ldots+\mu_{t-1}+1;n-m+\lambda_1+\lambda_2+\ldots+\lambda_{s-1}+1)\ne 0$.
       The only entries of matrix $A$ that are changed after the conjugation by $E_{t,1,j,1}(\xi)$, $t < j \leq k$,
       are the entries in the first row of the $j$-th block of rows. (Note that the column operation of conjugation $E_{t,1,j,1}(\xi)$
       does not change the matrix). Observe that $A_{11} \in \bkp{\underline{\mu}}{\underline{\mu}}$ and
       the only nonzero entries of $A_{12}$ lie in the top left corners of blocks according to partition $(\underline{\mu}, \underline{\lambda})$.
       Also, it is possible to choose $\xi$ such that
       $A(\mu_1+\mu_2+\ldots+\mu_{j-1}+1;n-m+\lambda_1+\lambda_2+\ldots+\lambda_{s-1}+1)=0$ for $j=t+1,t+2,\ldots,k$.

 \item[4.] Suppose that there exist $1 \leq t \leq k$ and $1 \leq s \leq l$ such that 
       $A(\mu_1+\mu_2+\ldots+\mu_{i-1}+1;n-m+\lambda_1+\lambda_2+\ldots+\lambda_{s-1}+1)=0$ for $i=1,2,\ldots,t-1$ and
       $A(\mu_1+\mu_2+\ldots+\mu_{t-1}+1;n-m+\lambda_1+\lambda_2+\ldots+\lambda_{s-1}+1)\ne 0$.
       Conjugation by $E_{k+j,1,k+s,1}(\xi)$, $s < j \leq l$,
       changes the first column in the $(k+j)$-th block of columns in the following way: 
       the entry $A(\mu_1+\mu_2+\ldots+\mu_{t-1}+1;n-m+\lambda_1+\lambda_2+\ldots+\lambda_{j-1}+1)$ vanishes for a suitable $\xi$ 
       and the entries above it do not change. However, the entries below it may change.
       At the same time conjugation by $E_{k+j,1,t,1}(\xi)$, $s < j \leq l$, makes (by the row operation) the entry 
       $A(n-m+\lambda_1+\lambda_2+\ldots+\lambda_{s-1}+1;n-m+\lambda_1+\lambda_2+\ldots+\lambda_{j-1}+2)\ne 0$.
       All entries outside of the $(k+j)$-th block of columns do not change.
       
       By the sequence of conjugations $E_{k+j,i,k+s,i}(\xi_i)$, $i=2,3, \ldots,\lambda_j$, all the entries in the $(k+s)$-th block
       of rows and $(k+j)$-th block of column become 0 and the only entries of the matrix that may change are the entries in 
       the last row of the $(k+s)$-th block of rows in $A_{21}$.
       
 \item[5.] If $A(\mu_1+\mu_2+\ldots+\mu_{t-1}+1;n-m+\lambda_1+\lambda_2+\ldots+\lambda_{s-1}+1)=0$, but 
       $A(\mu_1+\mu_2+\ldots+\mu_{j-1}+1;n-m+\lambda_1+\lambda_2+\ldots+\lambda_{s-1}+1)\ne 0$ for some $j>t$, then the conjugation
       by $E_{j,1,t,1}(1)$ causes that $A(\mu_1+\mu_2+\ldots+\mu_{t-1}+1;n-m+\lambda_1+\lambda_2+\ldots+\lambda_{s-1}+1)\ne 0$.
       The entries that are changed are all in $A_{11}$ and $A_{12}$, but $A_{11} \in \bkp{\underline{\mu}}{\underline{\mu}}$
       and the only nonzero entries of $A_{12}$ still lie in the top left corners 
       of blocks according to partition $(\underline{\mu}, \underline{\lambda})$.

 \item[6.] Suppose that $A(\mu_1+\mu_2+\ldots+\mu_{t-1}+1;n-m+\lambda_1+\lambda_2+\ldots+\lambda_{s-1}+1)\ne 0$ and let
       $E_{i,i}$ denote the matrix with the only nonzero entry $E_{i,i}(i;i)=1$. The conjugation by 
       matrix $I+\xi E_{\mu_1+\mu_2+\ldots+\mu_{t-1}+1,\mu_1+\mu_2+\ldots+\mu_{t-1}+1}$, for a suitable $\xi$, sets
       the entry $A(\mu_1+\mu_2+\ldots+\mu_{t-1}+1;n-m+\lambda_1+\lambda_2+\ldots+\lambda_{s-1}+1)$ to 1. 
       
 \item[7.] There exists a conjugation of $A$ by a suitable matrix such that it exchanges the
       $(k+i)$-th and $(k+j)$-th blocks of columns and $(k+i)$-th and $(k+j)$-th blocks of rows. Note that $A_{11}$ and $A_{22}$ 
       do not change. 
       
\end{enumerate}

From the description of conjugations 3., 4., 5. and 6. we see that they correspond to row and column operations on $A_{12}$ in the Gaussian 
elimination with a possible exchange of rows (but not columns). Therefore, it is possible 
to choose some sequence of conjugations, such that at the end the are exactly $\rk(A_{12})$ nonzero entries in $A_{12}$, lying
in first $\rk(A_{12})$ blocks of rows and they are all equal to 1. 

Denote the first column in each block of columns of submatrix $A_{12}$ by $X^{1}, X^{2},\ldots, X^{l}$.

For each $t$ and $s$  such that  $\lambda_{t-1}> \lambda_t=\lambda_{t+1}=\ldots=\lambda_{t+s-1}> \lambda_{t+s}$ and
$1 \leq t \leq t+s-1 \leq l$.
It is possible to choose a sequence of swaps, described in 7., such that the conjugated matrix has the property that there exists
an $i$, $0 \leq i \leq s-1$,  such that $\rk[X^{1}, X^2, \ldots, X^t]=\rk[X^{1}, X^2, \ldots, X^{t+i}]-i$ and 
$\rk[X^{1}, X^2, \ldots, X^{t+i}]=\rk[X^{1}, X^2, \ldots, X^{t+s-1}]$. I.e., the first $i$ columns of $X^{t}$, $X^{t+1}$,... $X^{t+s-1}$
are linearly independent and the other columns are equal to 0.

From now on, we will assume that matrix $A$ is already in the described form.

\bigskip

Similarly as in 3. and 5., it is possible to transform $A_{21}$ into its column echelon form without changing $A_{12}$ and $A_{22}$. 
Note that it would be also possible to make the row Gaussian elimination on $A_{21}$, but with these conjugations we would change the
form of $A_{12}$.

\bigskip

Using the above described conjugations, we transform the matrix $A=\left[ \begin{matrix}
 A_{11} & A_{12}\\
 A_{21} & J_{\underline{\lambda}}
\end{matrix} \right]$, so that it has the following properties:
\begin{enumerate}
 \item $A_{11} \in \sch{\underline{\mu}}{\underline{\mu}}$,
 \item $A_{12} \in \sch{\underline{\mu}}{1^m}$ and $A_{12}$ has
      $\rk\left[ X^{1}, X^{2}, \ldots,X^l\right]$ nonzero entries, all equal to 1. Moreover,
      $\rk[X^{1}, X^2, \ldots, X^t]=\rk[X^{1}, X^2, \ldots, X^{t+i}]-i$ and
      $\rk[X^{1}, X^2, \ldots, X^{t+i}]=\rk[X^{1}, X^2, \ldots, X^{t+s-1}]$ for all $t$ and $s$ such that
      $\lambda_{t-1}> \lambda_t=\lambda_{t+1}=\ldots=\lambda_{t+s-1}> \lambda_{t+s}$ and for an $i$, $0 \leq i \leq s$.
 \item $A_{21} \in \sch{1^m}{\underline{\mu}}$ is in a column echelon form with
       $A_{21}(i,j)\ne 0$  only if $i=n-m+\lambda_1+\lambda_2+\ldots+\lambda_t$ and $j=\mu_1+\mu_2+\ldots+\mu_s$
       for some $1 \leq t \leq l$ and $1 \leq s \leq k$.
\end{enumerate}
In this way we obtain a matrix $A$ that we shall call \DEF{a matrix in the reduced form}.

\bigskip 
 
\Example\label{ex:332}
Take $\underline{\mu}=(3,3,2)$ and $\underline{\lambda}=(3,2,2,1)$.
Then $A \in \sch{\underline{\mu}}{1^m}$,where $A_{22}=J_{(\lambda_1,\lambda_2,\ldots,\lambda_l)}$, has the following pattern:
$$ A= \left[ \begin{array}{ccc|ccc|cc||ccc|cc|cc|c}
 0 & 0 & x_1 & 0 & 0 & x_2 & 0 &  x_3 &  y_{11} & y_{12} & y_{13} & y_{21} & y_{22} & y_{31} & y_{32} & y_{41} \\
 0 & 0 & 0   & 0 & 0 & 0   & 0 &  0   &  0   & 0 & 0 & 0 & 0 & 0   & 0 & 0 \\
 0 & 0 & 0   & 0 & 0 & 0   & 0 &  0   &  0   & 0 & 0 & 0 & 0 & 0   & 0 & 0 \\
 \hline
 0 & 0 & x_4 & 0 & 0 & x_5 & 0 &  x_6 &  y_{51} & y_{52} & y_{53}  & y_{61} & y_{62} & y_{71} & y_{72} & y_{81} \\
 0 & 0 & 0   & 0 & 0 & 0   & 0 &  0   &  0   & 0 & 0 & 0 & 0 & 0   & 0 & 0 \\
 0 & 0 & 0   & 0 & 0 & 0   & 0 &  0   &  0   & 0 & 0 & 0 & 0 & 0   & 0 & 0 \\
 \hline
 0 & 0 & x_7 & 0 & 0 & x_8 & 0 &  x_9 &  y_{91} & y_{92} & y_{93} &y_{101}& y_{102} & y_{111}& y_{112} & y_{121} \\
 0 & 0 & 0   & 0 & 0 & 0   & 0 &  0   &  0   & 0 & 0 & 0 & 0 & 0   & 0 & 0 \\
 \hline
 0 & 0 & z_{11}  & 0 & 0 & z_{21} & 0 &  z_{31} &  0   & 1 & 0 & 0 & 0 & 0   & 0 & 0 \\
 0 & 0 & z_{12}  & 0 & 0 & z_{22} & 0 &  z_{32} &  0   & 0 & 1 & 0 & 0 & 0   & 0 & 0 \\
 0 & 0 & z_{13}  & 0 & 0 & z_{23} & 0 &  z_{33} &  0   & 0 & 0 & 0 & 0 & 0   & 0 & 0 \\
 \hline
 0 & 0 & z_{41} & 0 & 0 &  z_{51} & 0 & z_{61} &  0   & 0 & 0 & 0 & 1 & 0   & 0 & 0 \\
 0 & 0 & z_{42} & 0 & 0 &  z_{52} & 0 & z_{62} &  0   & 0 & 0 & 0 & 0 & 0   & 0 & 0 \\
 \hline
 0 & 0 & z_{71} & 0 & 0 & z_{81} & 0 &  z_{91} &  0   & 0 & 0 & 0 & 0 & 0   & 1 & 0 \\
 0 & 0 & z_{72} & 0 & 0 & z_{82} & 0 &  z_{92} &  0   & 0 & 0 & 0 & 0 & 0   & 0 & 0 \\
 \hline
 0 & 0 &z_{101}& 0 & 0 &z_{111}& 0 &z_{121}&  0   & 0 & 0 & 0 & 0 & 0   & 0 & 0 
\end{array}
\right]\, ,$$
where $x_i,y_j,z_{k} \in F$.

Suppose that
$\rk \left[ \begin{matrix}
y_{11} & y_{21}& y_{31}& y_{41}\\
y_{51} & y_{61}& y_{71}& y_{81}\\
y_{91} &y_{101}&y_{111}&y_{121}
\end{matrix} \right] =3$ and that $(y_{21},y_{61},y_{101})=\alpha(y_{11},y_{51},y_{91})$.

By the algorithm above, $A$ is first transformed to
\begin{equation*}
A=\left[ \begin{array}{ccc|ccc|cc||ccc|cc|cc|c}
 0 & 0 & x'_1 & 0 & 0 & x'_2 & 0 &  x'_3 &  1 & 0 & 0 & 0 & 0 & 0   & 0 & 0 \\
 0 & 0 & 0   & 0 & 0 & 0   & 0 &  0   &  0   & 0 & 0 & 0 & 0 & 0   & 0 & 0 \\
 0 & 0 & 0   & 0 & 0 & 0   & 0 &  0   &  0   & 0 & 0 & 0 & 0 & 0   & 0 & 0 \\
 \hline
 0 & 0 & x'_4 & 0 & 0 & x'_5 & 0 &  x'_6 &  0   & 0 & 0 & 1 & 0 & 0 & 0 & 0 \\
 0 & 0 & 0   & 0 & 0 & 0   & 0 &  0   &  0   & 0 & 0 & 0 & 0 & 0   & 0 & 0 \\
 0 & 0 & 0   & 0 & 0 & 0   & 0 &  0   &  0   & 0 & 0 & 0 & 0 & 0   & 0 & 0 \\
 \hline
 0 & 0 & x'_7 & 0 & 0 & x'_8 & 0 &  x'_9 &  0   & 0 & 0 & 0 & 0 & 0   & 0 & 1 \\
 0 & 0 & 0   & 0 & 0 & 0   & 0 &  0   &  0   & 0 & 0 & 0 & 0 & 0   & 0 & 0 \\
 \hline
 0 & 0 & 0   & 0 & 0 & 0   & 0 &  0   &  0   & 1 & 0 & 0 & 0 & 0   & 0 & 0 \\
 0 & 0 & 0   & 0 & 0 & 0   & 0 &  0   &  0   & 0 & 1 & 0 & 0 & 0   & 0 & 0 \\
 0 & 0 & z'_1 & 0 & 0 & z'_2 & 0 &  z'_3 &  0   & 0 & 0 & 0 & 0 & 0   & 0 & 0 \\
 \hline
 0 & 0 & 0   & 0 & 0 & 0   & 0 &  0   &  0   & 0 & 0 & 0 & 1 & 0   & 0 & 0 \\
 0 & 0 & z'_4 & 0 & 0 & z'_5 & 0 &  z'_6 &  0   & 0 & 0 & 0 & 0 & 0   & 0 & 0 \\
 \hline
 0 & 0 & 0   & 0 & 0 & 0   & 0 &  0   &  0   & 0 & 0 & 0 & 0 & 0   & 1 & 0 \\
 0 & 0 & z'_7 & 0 & 0 & z'_8 & 0 &  z'_9 &  0   & 0 & 0 & 0 & 0 & 0   & 0 & 0 \\
 \hline
 0 & 0 &z'_{10}& 0 & 0 &z'_{11}& 0 &z'_{12}&  0   & 0 & 0 & 0 & 0 & 0   & 0 & 0 
\end{array}
\right]\, . 
\end{equation*}
Depending on the values of $z'_i$'s, one of the possibilities for $A$ in the reduced form is for example
\begin{equation} \label{eq:a}
A=\left[ \begin{array}{ccc|ccc|cc||ccc|cc|cc|c}
 0 & 0 & x''_1 & 0 & 0 & x''_2 & 0 &  x''_3 &  1 & 0 & 0 & 0 & 0 & 0   & 0 & 0 \\
 0 & 0 & 0   & 0 & 0 & 0   & 0 &  0   &  0   & 0 & 0 & 0 & 0 & 0   & 0 & 0 \\
 0 & 0 & 0   & 0 & 0 & 0   & 0 &  0   &  0   & 0 & 0 & 0 & 0 & 0   & 0 & 0 \\
 \hline
 0 & 0 & x''_4 & 0 & 0 & x''_5 & 0 &  x''_6 &  0   & 0 & 0 & 1 & 0 & 0 & 0 & 0 \\
 0 & 0 & 0   & 0 & 0 & 0   & 0 &  0   &  0   & 0 & 0 & 0 & 0 & 0   & 0 & 0 \\
 0 & 0 & 0   & 0 & 0 & 0   & 0 &  0   &  0   & 0 & 0 & 0 & 0 & 0   & 0 & 0 \\
 \hline
 0 & 0 & x''_7 & 0 & 0 & x''_8 & 0 &  x''_9 &  0   & 0 & 0 & 0 & 0 & 0   & 0 & 1 \\
 0 & 0 & 0   & 0 & 0 & 0   & 0 &  0   &  0   & 0 & 0 & 0 & 0 & 0   & 0 & 0 \\
 \hline
 0 & 0 & 0   & 0 & 0 & 0   & 0 &  0   &  0   & 1 & 0 & 0 & 0 & 0   & 0 & 0 \\
 0 & 0 & 0   & 0 & 0 & 0   & 0 &  0   &  0   & 0 & 1 & 0 & 0 & 0   & 0 & 0 \\
 0 & 0 & z''_1 & 0 & 0 &  0 & 0 &  0 &  0   & 0 & 0 & 0 & 0 & 0   & 0 & 0 \\
 \hline
 0 & 0 & 0   & 0 & 0 & 0   & 0 &  0   &  0   & 0 & 0 & 0 & 1 & 0   & 0 & 0 \\
 0 & 0 &  z''_2 & 0 & 0 & 0 & 0 &  0 &  0   & 0 & 0 & 0 & 0 & 0   & 0 & 0 \\
 \hline
 0 & 0 & 0   & 0 & 0 & 0   & 0 &  0   &  0   & 0 & 0 & 0 & 0 & 0   & 1 & 0 \\
 0 & 0 & 0 & 0 & 0 & z''_3 & 0 &  0 &  0   & 0 & 0 & 0 & 0 & 0   & 0 & 0 \\
 \hline
 0 & 0 & 0 & 0 & 0 &  0 & 0 & z''_4&  0   & 0 & 0 & 0 & 0 & 0   & 0 & 0 
\end{array}
\right]\, ,
\end{equation}
where $z''_i \ne 0$ for $i=1,3,4$.
\hfill$\square$

\bigskip
\bigskip

\section{The Jordan canonical form of $A$} 

\bigskip

\Lemma
For a matrix $A\in \nbsch$ in the reduced form we have
$$A^{s+1}=\left[ \begin{matrix}
 A_{12}J_{\underline{\lambda}}^{s-1}A_{21} & A_{12}J_{\underline{\lambda}}^s\\
 J_{\underline{\lambda}}^s A_{21} & J_{\underline{\lambda}}^{s+1}
\end{matrix} \right]$$
for $s \geq 1$. 

\medskip

\begin{proof}
 By the properties of matrix $A$ it is clear that 
 $A_{11}^2=A_{11}A_{12}=A_{21}A_{11}=A_{21}A_{12}=0$.
 Now the expression for $A^{s+1}$ easily follows by induction on $s$.
\end{proof}

\bigskip

Note that $A_{12}$ has at most one nonzero entry (that is equal to 1) in every row, therefore all the entries of a matrix
$A^n$ are the entries of $A_{21}$. 

\bigskip

Recall that the \DEF{conjugated partition} of a partition $\underline{\lambda}$ is the partition 
$\underline{\lambda^T}=(\lambda_1^T, \lambda_2^T,\ldots,\lambda_{\lambda_1}^T)$,
where $\lambda_i^T=|\{j; \; \lambda_j \geq i \}|$ (and consequently $\lambda_i=|\{j; \; \lambda_j^T \geq i \}|$).  
It easily follows that $\lambda_{\lambda_i^T}=|\{j; \; \lambda_j^T \geq \lambda_i^T\}| \geq i$ and 
$\lambda_{\lambda_i^T+1}=|\{j; \; \lambda_j^T \geq \lambda_i^T+1\}| < i$. Therefore $\lambda_i^T$ is the maximal index $j$ such that
$\lambda_j \geq i$.

\smallskip

Recall that we denoted the first column in each block of columns of submatrix $A_{12}$ by $X^{1}, X^{2},\ldots, X^{l}$ and 
similarly denote the last rows in each block of rows of the submatrix $A_{21}$ by $Y_{1}, Y_{2},\ldots, Y_{l}$.
Moreover, denote by $Z_s^{1}, Z_s^{2},\ldots, Z_s^{k}$ the last columns in each block of columns of the submatrix 
$A_{12}J_{\underline{\lambda}}^{s-2}A_{21}$ of matrix $A^{s}$.

\smallskip

For $i=1,2, \ldots,l$ define $e_1(i)=\rk\left[X^1 X^2 \ldots X^{i} \right]$ 
and $e_2(i)=\rk \left[\begin{matrix} Y_1\\
Y_2\\
\vdots\\
Y_{i}
\end{matrix}\right]$.
From the patterns of $A_{12}$ and $A_{21}$, it follows that 
$\dim \im A_{12} J_{\underline{\lambda}}^s =e_1(\lambda_{s+1}^T)$ and
$\dim \im J_{\underline{\lambda}}^s A_{21}=e_2(\lambda_{s+1}^T)$.
Denote also 
$f(s)=\rk \left[ Z_{s}^{e_2(\lambda_{s}^T)+1} Z_{s}^{e_2(\lambda_{s}^T)+2} 
\ldots Z_{s}^{k}\right]$.
Now, it is clear that 
\begin{align*}
\rk\left(A^{s+1}\right)&=
\dim \im \left[\begin{matrix}
A_{12}J_{\underline{\lambda}}^s\\
J_{\underline{\lambda}}^{s+1}
 \end{matrix}\right]+
\dim \im \left[\begin{matrix}
 A_{12}J_{\underline{\lambda}}^{s-1}A_{21}\\
 J_{\underline{\lambda}}^s A_{21}
\end{matrix}\right]
=\\
&=\rk(J_{\underline{\lambda}}^{s+1})+\rk(A_{12}J_{\underline{\lambda}}^s)+\rk(J_{\underline{\lambda}}^sA_{12})+f(s+1)=\\
&=\lambda_{s+2}^T+\lambda_{s+3}^T+\ldots+\lambda_{\lambda_1}^T+e_1(\lambda_{s+1}^T)+e_2(\lambda_{s+1}^T)+f(s+1)\, .
\end{align*}

\medskip

\Lemma\label{thm:lengths}
The only possible lengths of Jordan chains of matrix $A$ are 1, 2, $\lambda_t$, $\lambda_t+1$ and $\lambda_t+2$ for some $t$.

\medskip

\begin{proof}
Without loss of generality suppose that $A$ is in its reduced form and
let $\{v_1,v_2,\ldots,v_n\}$ be the basis corresponding to matrix $A$. 
Denote $x(i)=v_{\mu_1+\mu_2+\ldots+\mu_{i}}$ for all $1 \leq i \leq k$. If $Ax(i) \notin \lin \{v_1,v_2,\ldots,v_{n-m}\}$, then
there exists $t(i)$, $1 \leq t(i) \leq l$, such that 
$Ax(i)=\sum\limits_{j=1}^k \alpha_j v_{\mu_1+\mu_2+\ldots+\mu_{j-1}+1}+\sum\limits_{j=t(i)}^l \beta_j v_{n-m+\lambda_1+\lambda_2+\ldots+\lambda_{j}}$
and $\beta_{t(i)} \ne 0$. 
It follows that $A^2x(i)=\sum\limits_{j\geq t(i), \lambda_j \geq 2} \beta_j v_{n-m+\lambda_1+\lambda_2+\ldots+\lambda_{j}-1}$,...,
$A^{\lambda_{t(i)}}x(i)=\sum\limits_{j\geq t(i), \lambda_j =\lambda_t(i)} \beta_j v_{n-m+\lambda_1+\lambda_2+\ldots+\lambda_{j-1}+1}$.

From the pattern of the submatrix $A_{12}$, it follows that $A^{\lambda_{t(i)}+1}x(i)=0$ if $e_1(t(i)-1)=e_1(t(i))$ and otherwise
$A^{\lambda_{t(i)}+1}x(i)=\sum\limits_{j=t(i)}^{t'(i)} \beta_j v_{\mu_1+\mu_2+\ldots+\mu_{j-s-1}+1}$  for some $0 \leq s \leq t(i)-1$.

Denote the set of all nonzero vectors $A^jx(i)$, $i=1,2,\ldots,k$, $j \geq 0$ by ${\cal A}$. 
The vectors in ${\cal A}$ are linearly independent and they form Jordan chains of length $\lambda_t(i)+2$, $\lambda_t(i)+1$, 2  and 1. 
Write $y(j)=v_{n-m+\lambda_1+\lambda_2+\ldots+\lambda_j}$ for $j \ne t(i)$, $1 \leq j \leq l$. Denote by ${\cal B}$ the set of all 
nonzero vectors of the form $y(j), Ay(j), A^2y(j),\ldots,A^{\lambda_{j}}y(j)$. It is clear that the vectors in ${\cal B}$ are 
linearly independent and form Jordan chains of length $\lambda_{j}+1$ and $\lambda_j$. 
From the pattern of the matrix $A$ it is easy to see that ${\cal A} \cup {\cal B}$ is a Jordan basis for $A$
and the only possible lengths of the Jordan chains of $A$ are 1, 2, $\lambda_t$, $\lambda_t+1$ and $\lambda_t+2$ for some $t$.
%
\end{proof}

\bigskip

We say that Jordan chains of matrix $A$ of lengths $\lambda_t$, $\lambda_t+1$ and $\lambda_t+2$ \DEF{arise from Jordan chain of 
$J_{\underline{\lambda}}$ of length $\lambda_t$}. By the proof of the Lemma \ref{thm:lengths} this is well defined.

\bigskip

\Lemma\label{thm:s+2}
Let $s=\lambda_t$ for some $1\leq t \leq l$. Then the number of Jordan chains of length $s+2$ 
of matrix $A$ that arise from Jordan chains of $J_{\underline{\lambda}}$ of length $s$
is equal to the $f(s+1)$.

\medskip

\begin{proof}
 For $s=\lambda_t$, $1 \leq t \leq l$, let $c_s$ be the number of Jordan chains of length $s+2$ that arise from Jordan chains of length $s$. 
 From Lemma \ref{thm:lengths} it follows that there exist indices $1 \leq t_1< t_2<  \ldots< t_{c_s}\leq l$ such that $\lambda_{t_i}=s$ and
 \begin{itemize}
  \item $A(n-m+\lambda_1+\lambda_2+\ldots+\lambda_{t_i}, \mu_1+\mu_2+\ldots\mu_{w_i}) \ne 0$ for some 
       $e_2(\lambda_{s+1}^T)+1 \leq w_1<w_2<\ldots<w_{c_s}\leq k$ and $i=1,2,\ldots,c_s$.
  \item $A(\mu_1+\mu_2+\ldots\mu_{w'_i-1}+1, n-m+\lambda_1+\lambda_2+\ldots+\lambda_{t_i-1}+1)\ne 0$ for some 
       $1 \leq w'_1<w'_2<\ldots<w'_{c_s}\leq k$ and $i=1,2,\ldots,c_s$.
 \end{itemize}
 Now, it follows that $(A_{12}J_{\underline{\lambda}}^{s-1}A_{21})(\mu_1+\mu_2+\ldots\mu_{w'_i-1}+1,\mu_1+\mu_2+\ldots\mu_{w_i}) \ne 0$
 for $i=1,2,\ldots,c_s$  and thus $f(s+1) \geq c_s$.
 
 If $f(s+1) \geq c_s+1$, then there exist
 $e_2(\lambda_{s+1}^T)+1 \leq w_0<w_1<w_2<\ldots<w_{c_s}\leq k$ and $1 \leq w'_0<w'_1<w'_2<\ldots<w'_{c_s}\leq k$, such that
 $(A_{12}J_{\underline{\lambda}}^{s-1}A_{21})(\mu_1+\mu_2+\ldots\mu_{w'_i-1}+1,\mu_1+\mu_2+\ldots\mu_{w_i}) \ne 0$
 for $i=0,1,2,\ldots,c_s$. Therefore, $A(n-m+\lambda_1+\lambda_2+\ldots+\lambda_{t_i}, \mu_1+\mu_2+\ldots\mu_{w_i}) \ne 0$ for
 and $A(\mu_1+\mu_2+\ldots\mu_{w'_i-1}+1, n-m+\lambda_1+\lambda_2+\ldots+\lambda_{t_i-1}+1)\ne 0$ for $i=0,1,2,\ldots,c_s$ and
 the number of Jordan chains of length $s+2$ that arise from Jordan chains of length $s$ would be at least $c_s+1$. This contradicts the
 definition of the number $c_s$ and thus the Lemma follows.
\end{proof}

\bigskip

\Lemma\label{thm:s+1}
Let $s=\lambda_t$ for some $1\leq t \leq l$. Then the number of Jordan chains of length $s+1$ 
of matrix $A$ that arise from Jordan chains of $J_{\underline{\lambda}}$ of length $s$ is equal to 
$e_1(\lambda_s^T)-e_1(\lambda_{s+1}^T)+e_2(\lambda_s^T)-e_2(\lambda_{s+1}^T)-2f(s+1)$.

\medskip

\begin{proof}
 From  Lemmas \ref{thm:lengths} and \ref{thm:s+2}  it follows that the number of 
 Jordan chains of length $s+1$, that arise from a Jordan chain of length $s$ and a nonzero element in $A_{21}$, is equal to 
 $e_2(\lambda_s^T)-e_2(\lambda_{s+1}^T)-f(s+1)$. Similarly, the number of Jordan chains of length $s+1$, that arise from a Jordan chain of length 
 $s$ and a nonzero element in $A_{12}$, is equal to $e_1(\lambda_s^T)-e_1(\lambda_{s+1}^T)-f(s+1)$. 
 Again, by Lemma \ref{thm:lengths}, there are no other Jordan 
 chains of length $s+1$ and this proves the Lemma.
\end{proof}

\bigskip

\Example
Let us recall the case $\sh{B}=(3,3,2,1^8)$ and $\underline{\lambda}=(3,2,2,1)$ as in Example \ref{ex:332} and take a matrix
$A \in \nbsch$ in the reduced form as in \eqref{eq:a}. 
Write $g(s+1)=e_1(\lambda_s^T)-e_1(\lambda_{s+1}^T)+e_2(\lambda_s^T)-e_2(\lambda_{s+1}^T)-2f(s+1)$.
It can be easily computed that $f(2)=1$, $g(2)=0$, $f(3)=0$, $g(3)=2$, $f(4)=1$ and $g(4)=0$. Therefore $\sh{A}=(5,3^3,1^2)$.

\bigskip

\Theorem\label{thm:nbsch}
If $\sh{B}=(\mu_1,\mu_2,\ldots,\mu_k,1^m) \in \partition{n}$, then
\begin{align*}
 \partition{\nbsch}=&\left\{ \ord{ \lambda_1+\varepsilon_1,\lambda_2+\varepsilon_2,\ldots,\lambda_l+\varepsilon_l, 
         2^{c}, 1^{d}} ; \right. \; \underline{\lambda}\in \partition{m}, \varepsilon_i \in \{0,1,2\},  \\
         & \left. \quad  0 \leq 2c \leq  2k-\sum_{i=1}^l \varepsilon_i, 
         2c+d+\sum_{i=1}^l \left( \lambda_i+\varepsilon_i\right)=n
\right\} \, .
\end{align*}

\medskip

\begin{proof} 
  For an arbitrary matrix $A \in \nbsch$ it is clear from Lemmas \ref{thm:lengths}, \ref{thm:s+2} and \ref{thm:s+1} that 
  $\sh{A}=\left(\lambda_1+\varepsilon_1,\lambda_2+\varepsilon_2,\ldots,\lambda_l+\varepsilon_l, 
         2^{c}, 1^{d}\right)$, where $\varepsilon_i \in \{0,1,2\}$. 
  By Lemmas \ref{thm:s+2} and \ref{thm:s+1},
  $\sum_{i=1}^l \varepsilon_i =e_1(\lambda_1^T)+e_2(\lambda_1^T)-\sum_{s=1}^{\lambda_l}f(s+1)$. 
  Denote by $A_{11}^{i}$ the $i$-th column of submatrix $A_{11}$. Since 
  $c=
   \left[ 
   A_{11}^{\mu_1+\mu_2+\ldots+\mu_{e_2(l)}+1}, A_{11}^{\mu_1+\mu_2+\ldots+\mu_{e_2(l)}+2}, \ldots,A_{11}^{\mu_1+\mu_2+\ldots+\mu_{k}} 
   \right]$,
  it follows that $c+e_2(l)\leq k$. Similarly, $c+e_1(l) \leq k$, and thus 
  $2k-2c \geq e_1(l)+e_2(l) \geq e_1(\lambda_1^T)+e_2(\lambda_1^T)-\sum_{s=1}^{\lambda_l}f(s+1)=\sum_{i=1}^l \varepsilon_i$.
  
  Conversely, let $\underline{\nu}=\ord{ \lambda_1+\varepsilon_1,\lambda_2+\varepsilon_2,\ldots,\lambda_l+\varepsilon_l, 
         2^{c}, 1^{d}} \in \partition{n}$, such that
  $\underline{\lambda} \in \partition{m}$, $\varepsilon_i \in \{0,1,2\}$ 
  and $0 \leq c \leq k- \frac{1}{2} \sum_{i=1}^l \varepsilon_i$. We will construct a matrix $A$, such that 
  $\sh{A}=\underline{\nu}$.
  
  
  Let us denote $\varphi_2(i)=|\{j > i; \; \varepsilon_j=2\}|$. 
  Let $t_0=s_0=0$ and for $i=1,2, \ldots,l$ let 
  $$t_i=\left\{ \begin{array}{cl}
      t_{i-1}+1, &  \text{if } \varepsilon_i=1  \text{ and }  t_{i-1}+1+\varphi_2(i) \leq k-c \text{ or } \varepsilon_i=2,\\
      t_{i-1}, & \text{otherwise }.
    \end{array}
  \right.$$
  and
  $$s_i=\left\{ \begin{array}{cl}
    s_{i-1}+1, &  \text{if } \varepsilon_i = 1 \text{ and } t_{i}=t_{i-1} \text{ or } \varepsilon_i=2,\\
    s_{i-1},& \text{otherwise }.
  \end{array}
  \right.$$
  Now, define matrix $A \in \nbsch$ with the following properties.
  \begin{itemize}
   \item If $t_{i-1}<t_{i}$ for some $1 \leq i < l$, then let $A(n-m+\lambda_1+\lambda_2+\ldots+\lambda_i,\mu_1+\mu_2+\ldots+\mu_{t_i})=1$.
   \item If $s_{i-1}<s_{i}$ for some $1 \leq i < l$, then let $A(\mu_1+\mu_2+\ldots+\mu_{s_i-1}+1,n-m+\lambda_1+\lambda_2+\ldots+\lambda_{i-1}+1)=1$.
   \item $A_{22}=J_{\lambda_1}\oplus J_{\lambda_2}\oplus \ldots \oplus J_{\lambda_l}$.
   \item For $i=1,2,\ldots,a$, let $A(\mu_1+\mu_2+\ldots+\mu_{k-i}+1,\mu_1+\mu_2+\ldots+\mu_{k-i+1})=1$.
  \end{itemize}
  Now, it can be easily seen that $\sh{A}=\ord{ \lambda_1+\varepsilon_1,\lambda_2+\varepsilon_2,\ldots,\lambda_l+\varepsilon_l, 2^{c}, 1^{d}}$.
\end{proof}
  
\bigskip

Recall that $\nschab$ is the variety of all pairs 
$(A,B)\in \nsch$ such that $AB=BA=A^a=B^b=0$ and that $\pi_1, \pi_2: \nschab \to \n$ are the projection maps, i.e. 
$\pi_1(A,B)=A$ and $\pi_2(A,B)=B$ for $(A,B) \in \nschab$.

Then $\partition{\nschab}$ is equal to the set of all pairs
of partitions $\underline{\mu}, \underline{\nu} \in \partition{n}$, such that $\mu_1 \leq a$, $\nu_1 \leq b$
and the intersection of fibers $\pi_1^{-1}({\cal O}_{\underline{\mu}}) \cap \pi_2^{-1}({\cal O}_{\underline{\nu}})$
is nonempty.
The following is an easy consequence of Theorem \ref{thm:nbsch}.

\bigskip

\Corollary
\begin{align*}
 \partition{\nschab}=&\bigg\{  (\underline{\mu}, \underline{\nu})\in \partition{n} \times \partition{n}; \;
    \underline{\mu}=(\mu_1,\mu_2,\ldots,\mu_k,1^m),  \, \mu_1 \leq a, \, \mu_k \geq 2,  \\ 
   & \qquad \underline{\nu}=\ord{ \lambda_1+\varepsilon_1,\lambda_2+\varepsilon_2,\ldots,\lambda_l+\varepsilon_l, 2^{c}, 1^{d}}, \,
      \underline{\lambda} \in \partition{m}, \, \varepsilon_i \in \{0,1,2\},   \\
   & \qquad  \lambda_i+\varepsilon_i \leq b, \, 0 \leq 2c \leq  2k-\sum_{i=1}^l \varepsilon_i, \; 
      2c+d+\sum_{i=1}^l \left( \lambda_i+\varepsilon_i\right)=n  \bigg\} 
   \tag*{$\blacksquare$}
\end{align*}

\bigskip

It was proved by Schr\"oer \cite[Theorem 1.1]{schroeer} that the irreducible components of $\nsch$ are
$${\cal C}_j=\{(A,B) \in \nsch ; \; \rk(A) \leq n-j, \, \rk(B) \leq j\}$$
for $1 \leq j \leq n-1$. Using Theorem \ref{thm:nbsch} it is now easy to characterize the set $\partition{{\cal C}_j}$.

\bigskip

\Corollary
\begin{align*}
 \partition{{\cal C}_j}=
 &\bigg\{ (\underline{\mu}, \underline{\nu})\in \partition{n} \times \partition{n}; \;
    \underline{\mu}=(\mu_1,\mu_2,\ldots,\mu_k,1^m), \, \mu_k \geq 2, \, \underline{\lambda} \in \partition{m},  \\ 
   & \qquad \underline{\nu}=\ord{ \lambda_1+\varepsilon_1,\lambda_2+\varepsilon_2,\ldots,\lambda_l+\varepsilon_l, 2^{c}, 1^{d}}, \,
       \, \varepsilon_i \in \{0,1,2\},   \\
   & \qquad 0 \leq 2c \leq  2k-\sum_{i=1}^l \varepsilon_i, \; 
      2c+d+\sum_{i=1}^l \left( \lambda_i+\varepsilon_i\right)=n, \\
   & \qquad  n-l-c-d \leq j \leq k+m \bigg\} 
   \tag*{$\blacksquare$}
\end{align*}

\end{document}